\documentclass[11pt]{article}

\usepackage{graphics}
\usepackage{graphicx}
\usepackage{amssymb}
\usepackage{amsmath}
\usepackage{epsf}
\usepackage{verbatim}

\newcommand{\halmos}{\rule{1ex}{1.4ex}}
\newcommand{\qed}{\hfill \halmos} %put !#at right margin
\newcommand{\be}[1]{\begin{equation}\label{#1}}
\newcommand{\ee}{\end{equation}}
\newcommand{\bi}{\begin{itemize}}
\newcommand{\ei}{\end{itemize}}
\newcommand{\ben}{\begin{enumerate}}
\newcommand{\een}{\end{enumerate}}

\newcommand{\R}{{\mathbb R}}  %ams bold

\newcommand{\bl}[1]{\begin{Lemma}\label{#1}}
\newcommand{\el}{\qed\end{Lemma}}
\newcommand{\bt}[1]{\begin{Theorem}\label{#1}}
\newcommand{\et}{\end{Theorem}}
\newcommand{\epr}{\end{proof}}
\newcommand{\bpr}{\begin{proof}}
\newenvironment{proof}{\noindent {\em Proof}.\ }{\hspace*{\fill}$\halmos$\medskip}

\newcommand{\beqn}{\begin{eqnarray*}}
\newcommand{\eeqn}{\end{eqnarray*}}

\newtheorem{Theorem}{Theorem}
\newtheorem{Proposition}[Theorem]{Proposition}
\newtheorem{Lemma}[Theorem]{Lemma}
\newtheorem{Claim}[Theorem]{Claim}
\newtheorem{Corollary}[Theorem]{Corollary}

\newtheorem{Question}[Theorem]{Question}

\newtheorem{Example}[Theorem]{Example}

\newcommand{\norma}[1]{\ensuremath{\left| #1 \right|}}

\newcommand{\rest}{\upharpoonright}

\newcommand{\flo}[1]{\left\lfloor #1 \right\rfloor}

\begin{document}

\title{Analogues of the Smale and Hirsch Theorems for Cooperative Boolean and Other Discrete Systems}

\author{
\and Winfried
Just\footnote{Department of Mathematics, Ohio University} \ and
G.A. Enciso\footnote{Mathematical Biosciences Institute, Ohio State University, and Harvard Medical School, Department of Systems Biology.
This material is based upon work supported by the National Science Foundation
under Agreement No. 0112050 and by The Ohio State University.}
}

\maketitle

\begin{center}
 \emph{Dedicated to Avner Friedman, on the occasion of his 75th birthday.}
\end{center}

\begin{abstract}
Discrete dynamical systems defined on the state space $\Pi=\{0,1,\ldots,
p-1\}^n$ have been used in multiple applications, most recently for the
modeling of gene and protein networks.   In this paper we study to what extent
well-known theorems by Smale and Hirsch, which form part of the theory of
(continuous) monotone dynamical systems, generalize or fail to do so in the
discrete case.

We show that that arbitrary $m$-dimensional systems cannot necessarily be
embedded into $n$-dimensional cooperative systems for $n=m+1$, as in the Smale
theorem for the continuous case, but we show that this is possible for $n=m+2$
as long as $p$ is sufficiently large.

We also prove that a natural discrete analogue of strong cooperativity implies nontrivial
bounds on the lengths of periodic orbits and imposes a condition akin to Lyapunov
stability on all attractors.  Finally, we explore several natural candidates for definitions of
irreducibility of a discrete system.
While some of these notions imply the strong cooperativity of a given cooperative system
 and impose even tighter bounds on the lengths of periodic orbits
than strong cooperativity alone, other plausible definitions allow the existence of exponentially
long periodic orbits.

\end{abstract}

%Vinny: add keywords
%\paragraph{Keywords:} These are the keywords


\vspace{1ex}


%92B99  mathematical biology, other
%39A11 Stability and asymptotics of difference equations; oscillatory and periodic solutions
%34C12 Monotone systems

\bigskip

%\vspace{1in}
%gae4 why the space??

\section{Introduction}

Let $(L, \leq)$  be a linearly ordered set, let $n \geq 1$, let $L_1, \ldots L_n \subseteq L$ with the induced order,
and consider the set $\Pi= \prod_{i=1}^n L_i$.
A map $g:\Pi\to \Pi$ defines the discrete dynamical system

\be{discrete} x(t+1)=g(x(t)), \ \ \ x(t)\in \Pi, \ee

We call (\ref{discrete}) an \emph{$n$-dimensional, discrete system} and
also identify it with the pair $(\Pi,g)$.  For most of this paper, $L$ will be the set of real numbers with the natural order,
and $L_i = \{0, ... , p-1\}$ for some fixed integer $p > 1$.  In this case we speak of an \emph{$n$-dimensional, $p$-discrete system.}
The case $p=2$ corresponds to the
so-called \emph{Boolean networks} or \emph{Boolean systems} which are used in various  disciplines,
notably in the study of gene regulatory systems
\cite{Kauffman:N1969,Kauffman:JTB1969,Tang:PhRE2007,Albert:JTB2003,
Samal:Jain:A2007, Shmulevich:IEEE2002,Sontag:SSB2007, Zaliapin:JSP2004}.
If all $L_i$'s are finite, then we may without loss of generality assume that
$L_i = \{0, \ldots , p_i-1\}$ for some $p_i > 1$, but the $p_i$'s are not necessarily all equal.  In this case we speak of a
\emph{finite discrete system.}

 Define a partial order on $\Pi$ by $x\leq y$ if $x_i\leq y_i$ for $i=1,\ldots, n$.  We call this relation the \emph{cooperative order,} and we will
 not make a notational distinction between it and the order relation on $L$.
A discrete system (\ref{discrete}) is said to be \emph{cooperative} if $x(0)\leq y(0)$ implies  $x(t)\leq y(t)$ for every $t\geq0$, where $x(t),y(t)$ are the solutions of the system with initial conditions $x(0),y(0)$ respectively.
Clearly this is equivalent to the property that $x\leq y$ implies $g(x)\leq
g(y)$.
 Discrete cooperative systems have been proposed as a tool to study genetic networks by Sontag and others \cite{Sontag:Laubenbacher:A2007,Sontag:SSB2007}.

The cooperativity property has a well-studied counterpart in continuous
dynamical systems \be{ode} \frac{d x_i}{d t}=f_i(x),\ \ \ i=1\ldots n,
 \ee
for $C^1$ vector fields $f:\R^n\to \R^n$.  Namely, the system (\ref{ode}) is
\emph{cooperative} if whenever $x(t),\ y(t)$ are two solutions such that
$x_i(0)\leq y_i(0)$, $i=1,\ldots, n$, then $x_i(t)\leq y_i(t)$ for every $t>0$,
$i=1,\ldots, n$.  Cooperative systems are canonical examples of so-called
monotone systems, which have been studied extensively by M. Hirsch, H. Smith,
H. Matano, P. Pol\'{a}\v{c}ik and others, and more recently by Sontag and
collaborators in the context of gene regulatory networks under exclusively
positive feedback \cite
{Sontag:multi,Hirsch:1985,Smith:review,Polacik:handbook,Smith:monotone}.

\subsection*{The Smale and Hirsch Theorems}

In the present paper we consider two important results from the theory of
(continuous) monotone dynamical systems, and we show to what extent these
results either generalize or fail to do so in the context of cooperative
discrete systems (\ref{discrete}).

The first result was originally published by S. Smale in the 1970's
\cite{Smale:1976}. It states in this context
that any compactly supported, $(n-1)$-dimensional, $C^1$ dynamical system
defined on $H=\{x\in \R^n\,|\, x_1+\ldots + x_n=0\}$ can be embedded into some
cooperative $C^1$ system~(\ref{ode}).  Equivalently, the dynamics of
cooperative systems can be completely arbitrary on unordered hyperplanes such
as $H$.  See also \cite{Enciso:continuum}, where the cooperative system
(\ref{ode}) is shown to have bounded solutions and only two equilibria outside
of $H$.

One way to regard the Smale theorem in the discrete case would be to ask
whether discrete cooperative systems can have arbitrary dynamics on unordered
sets $H=\{x\in \Pi\,|\, x_1+\ldots + x_n=\mbox{const.}\}$.   This is trivially
true; see Lemma~\ref{lemma trivial embedding}.

 An alternative approach is to study whether one can embed an arbitrary $m$-dimensional $p$-discrete system (\ref{discrete}) into a cooperative
 $(m+1)$-dimensional $p$-discrete system.  We show that the answer to this question is \emph{no} (Theorem~\ref{theorem bounds}, item~3), but that the statement is true (for sufficiently large $p$) if `$m+1$' is weakened to `$m+2$' (Theorem~\ref{theorem bounds}, item~2).

The second result for continuous cooperative systems was proved by M. Hirsch in
\cite{Hirsch:1985}.  A continuous cooperative system is \emph{strongly
cooperative} if for every two different initial conditions $x(0) \leq y(0)$ we
have $x_i(t) < y_i(t)$ for all $i = 1, \ldots , n$ and $t > 0$.
A closely related definition involves the digraph $G$ associated with the
system:  in the cooperative case, $G$ is defined as having nodes $1,\ldots, n$,
and an arc from $i$ to $j$ is present if and only if $\partial f_j /\partial
x_i(x)>0$ on $\R^n$.  A continuous cooperative system (\ref{ode}) is strongly
cooperative if the digraph $G$ is strongly connected \cite{Smith:monotone}; we
refer to the latter property as the \emph{irreducibility} of the system
(\ref{ode}). Hirsch's theorem states that almost every bounded solution of a
strongly cooperative system (\ref{ode}) converges towards the set of
equilibria. This result rules out stable periodic orbits and chaotic behavior.
It was also generalized for abstract order relations in
Banach spaces by Hirsch and extended to continuous-space, discrete-time maps by
Tere\v{s}\v{c}\'{a}k, Pol\'{a}\v{c}ik and collaborators; see
\cite{Hirsch:1988,Polacik:handbook,Terescak:1992,
Enciso:Hirsch:Smith:JDDE2006}.

For finite discrete systems, we will consider analog definitions of strong
cooperativity and of irreducibility of a cooperative system (\ref{discrete}).
We are particularly interested in whether these definitions rule out the existence of
exponentially long periodic orbits,
which in finite discrete systems can be considered analogues of
chaotic attractors.
  We show that strong cooperativity does not rule out
periodic orbits altogether, but that it puts a nontrivial, subexponential bound
on their lengths and imposes a condition akin to Lyapunov
stability on all attractors.  Finally, we explore several natural candidates for definitions of
irreducibility of a finite discrete system.  We show that predicted properties of the system can dramatically change when
subtle changes to our definitions are made.  While some plausible definitions of irreducibility still allow for exponentially
long periodic orbits (and hence do not imply strong cooperativity), other definitions of
irreducible cooperative systems imply strong cooperativity and impose a bound of $n$ (the dimension of the system)
on the lengths of periodic orbits.  This is a much tighter bound than the one implied by strong
cooperativity alone.

\subsection*{Outline of the Sections}

In Section~\ref{section unordered} we give a general condition under which
an arbitrary $m$-dimensional, $p$-discrete system can be embedded into a
cooperative $n$-dimensional, $q$-discrete system (Proposition~\ref{proposition
embedding}).  We rely on several standard results from the literature,
especially a generalization of the classical Sperner theorem.  In
Section~\ref{section bounds} we provide bounds on the maximum size $d_{n,p}$ of
an unordered subset of $\Pi$, and we use these bounds to study the special
cases $n=m+1$ and $n=m+2$ (Theorem~\ref{theorem bounds}).  In
Section~\ref{smaleapprox}  we prove a general result on extensions of
cooperative partial functions on $\Pi$ to cooperative systems on
$\Pi$
and discuss how our results are related to a certain
generalization of Smale's theorem.
We give a short discussion in Section~\ref{almost coop section} about applying Theorem~\ref{theorem
bounds} to the case of \emph{almost cooperative} discrete systems
\cite{Sontag:SSB2007}, by showing a simple example of an almost
cooperative Boolean system of dimension $m$ that cannot be embedded into a
cooperative Boolean system of dimension $m+1$.
In Section~\ref{section alternative} we introduce a counterpart
of  strong cooperativity
 for finite discrete systems (\ref{discrete}) and show that
it imposes substantial restrictions on the possible dynamics.  In particular, we show that
strongly cooperative $p$-discrete systems cannot have exponentially long periodic orbits.
 In Section~\ref{section strong coop} we explore several natural definitions of
 irreducibility for finite discrete systems and prove bounds on the lengths of periodic orbits in cooperative systems
 that are irreducible in the sense of these definitions.

\section{Unordered Sets and Cooperative Embeddings of $p$-discrete systems}   \label{section unordered}

 Let $\Sigma:=\prod_{i=1}^m L_i$ and $\Pi = \prod_{i=1}^n L^*_i$ and consider an arbitrary map $f:\Sigma\to \Sigma$.  A \emph{cooperative embedding} of $(\Sigma,f)$ into a cooperative system $(\Pi,g)$ as in~(\ref{discrete}) is an injective function $\phi:\Sigma\to \Pi$ such that $g(\phi(x))=\phi(f(x))$ for every $x\in \Sigma$.
  If $\Pi = \prod_{i=1}^n \{0, \ldots ,p_i -1\}$, then we define $S(x)=x_1+\ldots + x_n$ for $x\in
\Pi$.  These definitions will be used throughout this paper.

For the remainder of this section and the next one, let $\Pi = \{0, \ldots , p-1\}^n$
%gae4 replaced p with p-1 above
 for some fixed integer $p > 1$.
We will compute the least dimension $n$ such that any $m$-dimensional $p$-discrete system $(\Sigma, f)$ can be embedded into an $n$-dimensional cooperative system $(\Pi, g)$.

   A subset $A\subseteq
\Pi$ is said to be \emph{unordered} if no two different elements $a,b\in A$
satisfy $a\leq b$.  Define the set

\begin{equation}\label{Ddef}
D:=\{x\in \Pi\, |\, S(x) = \flo{ n (p-1) /2}\}, \ \ \ d_{n,p}:=\norma{D}.
\end{equation}
%
%gae4 wrote "Define the set" above
This set $D$ is clearly unordered, because if $x\leq y$, $x\not=y$, then
necessarily  $S(x)<S(y)$, and $x,y$ cannot  be both in $D$. Notice that
$d_{n,2}=\binom{n}{\flo{n/2}}$.   We quote a generalization of Sperner's
Theorem  \cite{Anderson:2002,Clements:1984},
 which states that $D$ is a set of maximum size in $\Pi$ with this property:

\begin{Lemma}  \label{lemma Sperner}
Consider the set $\Pi=\{0,1,\ldots,p-1\}^n$, under the cooperative order
$\leq$. Then $|A| \leq d_{n,p}$, for any unordered set $A$.
\end{Lemma}

The following lemma will be important below, see Proposition~5.2 in
\cite{Smith:review} for a proof.

\begin{Lemma} \label{lemma unordered}
Consider a cooperative map $g$ defined on a space $\Pi$.  Then any periodic
orbit is unordered.
\end{Lemma}

Another basic property of unordered sets is the following `trivial embedding'
result, which is well known at least for the Boolean case (see for instance
\cite{Sontag:SSB2007}).

\begin{Lemma}  \label{lemma trivial embedding}
Let $A\subseteq \Pi$ be unordered, and let $\gamma:A\to A$ be an arbitrary
function.  Then there exists a cooperative system (\ref{discrete}) such that
$g|_A=\gamma$.
\end{Lemma}

\bpr Let $\hat{A}$ be any unordered subset of $\Pi$ which contains $A$, and
which is maximal with respect to this property.  Define $g(a):=\gamma(a)$ for
$a\in A$, and $g(a)=a$ for $a\in \hat{A}-A$.  For all other $x\in \Pi$, there
must exist $a\in\hat{A}$ such that either $a\leq x$ or $x\leq a$, by the
maximality of $\hat{A}$.  If $x\leq a$ let $g(x):=[0,\ldots,0]$, and if $a\leq
x$ let $g(x):=[p-1,\ldots,p-1]$. \epr

\begin{Proposition}\label{proposition embedding}
Let $n,m$ be positive integers, $p,q>1$, and $\Pi=\{0,1,\ldots, p-1\}^n$,
$\Sigma=\{0,1,\ldots, q-1\}^m$.  Then the following are equivalent:
\begin{description}
\item{(i)} Any discrete system $(\Sigma, f)$ can be embedded into a cooperative discrete system $(\Pi, g)$.

\item{(ii)} $q^m\leq d_{n,p}$.
\end{description}
\end{Proposition}

\bpr Suppose first that $q^m\leq d_{n,p}$, and consider any discrete system
$(\Sigma, f)$.   We use an arbitrary injective function $\phi:\Sigma\to \Pi$
such that $A:=\mbox{Im}(\phi)\subseteq D$.  Let $\gamma(y):=\phi(f(x))$
whenever $y\in A$, where $x=\phi^{-1}(y)$. Thus by construction
$\gamma(\phi(x))=\gamma(y)=\phi(f(x))$ holds for $x\in \Sigma$.    Apply
Lemma~\ref{lemma trivial embedding} to define $g$ and obtain a full cooperative
embedding.

Now assume (i) in the statement.  To prove that (ii) must hold, simply consider
a map~$f$ on $\Sigma$ which generates a single orbit with period $q^m$.  By
(i), there exists an embedding into~$\Pi$, and the image of $\Sigma$ is
unordered in $\Pi$ by Lemma~\ref{lemma unordered}.  The inequality follows from
Lemma~\ref{lemma Sperner}. \epr

Another form of cooperative embedding was given by Smith
\cite{Smith:discreteSGT} for a large class of non-cooperative, but possibly
continuous maps.   In that case $n=2m$ holds.  By Proposition~\ref{proposition
embedding}, a much sharper bound holds for the discrete case.

\section{Bounds on Discrete Cooperative Embeddings}  \label{section bounds}

Let $p>1$, $n>0$ be arbitrary, and let $\Pi,D,d_{n,p}$ be as in the previous section.  We begin
this section with several lemmas.

\begin{Lemma}\label{lemma dnp lower bound}
$\displaystyle d_{n,p}\geq \frac{p^{n-1}}{n}$.
\end{Lemma}

\bpr Let $S_j:=\{x\in \Pi\,|\, S(x)=j\}$, for $j=0,\ldots, n(p-1)$.  Each of
these sets is unordered, and therefore $\norma{S_j}\leq d_{n,p}$ by
Lemma~\ref{lemma Sperner}. Therefore
\[
p^n = \sum_{j=0}^{n(p-1)} \norma{S_j} \leq (n(p-1)+1)d_{n,p}\leq n p d_{n,p}.
\]
\epr

\begin{Lemma}\label{expolemma}
 Let $c$ be such that  $0<c<p$.  Then  $d_{n,p}\geq c^n$, for all sufficiently large $n$.
\end{Lemma}

\bpr By Lemma~\ref{lemma dnp lower bound}, it is sufficient to show that
$p^{n-1}/n \geq c^n$.  But this is equivalent to $\ln p \geq \ln c + (\ln n +
\ln p)/n$.  This inequality holds for large $n$ since $\ln p > \ln c$. \epr

We now prove an upper bound for $d_{n,p}$.

\begin{Lemma}\label{dnpupper}
Let $p,n > 1$. Then $d_{n+1, p} < p^n$.
\end{Lemma}

\bpr
Let $p, n$ be as in the assumptions, and let $x$ be a randomly chosen element of
$\{0, \ldots , p-1\}^{n+1}$ with the uniform distribution.  For $x$ to be in $D$, we must have
$\lfloor n(p-1)/2\rfloor - p + 1 \leq x_1 + \dots + x_{n-1} \leq \lfloor n(p-1)/2\rfloor$ and
$x_n =  x_1 + \dots + x_{n-1} - \lfloor n(p-1)/2\rfloor$.  Let
$A$ be the event that $\lfloor n(p-1)/2\rfloor - p + 1 \leq x_1 + \dots + x_{n-1} \leq \lfloor n(p-1)/2\rfloor$.
Our assumption on $n$ implies that $P(A) < 1$.  Moreover, note that
$P(x_n =  x_1 + \dots + x_{n-1} - \lfloor n(p-1)/2\rfloor|A) = \frac{1}{p}$.
Thus $P(x \in D) = \frac{d_{n+1, p}}{p^{n+1}} < \frac{1}{p}$, and the lemma follows.
\epr

The above estimates have important consequences for embeddings of
$m$-dimensional finite discrete systems into $n$-dimensional cooperative finite discrete systems.  In particular, unlike for continuous systems, for large $m$, an
$m$-dimensional $p$-discrete system can in general not be embedded into an
$(m+1)$-dimensional $p$-discrete cooperative system.

\begin{Theorem}  \label{theorem bounds}
The following statements hold:

\begin{enumerate}
\item For every $p>1$, and for every $m>0$, there exists $n>m$ such that every $m$-dimensional $p$-discrete system
can be embedded into an $n$-dimensional cooperative $p$-discrete system.

\item For every $m>0$, there exists $p_0$ such that for every $p > p_0$  every $m$-dimensional $p$-discrete system can be embedded into a cooperative $p$-discrete system of dimension $m+2$.

\item For every $m, p > 1$  there exists an $m$-dimensional $p$-discrete system  that cannot be embedded into a cooperative $p$-discrete system of dimension $m+1$.
\end{enumerate}
\end{Theorem}

\bpr The first two statements are immediate consequences of Lemma~\ref{lemma
dnp lower bound} and Proposition~\ref{proposition embedding}.  For the first
one, let $n$ be large enough so that $p^m\leq p^{n-1}/n$.  Then $p^m\leq d_{n,p}$,
and the conclusion follows.  For the second statement, let simply $p\geq m+2$.
Then
\[
p^m\leq \frac{p^{m+1}}{m+2}=\frac{p^{(m+2)-1}}{m+2}\leq d_{m+2,p}.
\]
For the third statement, let $m,p > 1$.
Let $f$ be defined on $\Pi$ so as to generate a single orbit of length $p^m$.
Then the image of $\Pi$ under any embedding $\phi$ into
$\Sigma=\{0,\ldots,p-1\}^{m+1}$ would also generate a periodic orbit of this
length.  Assuming that the system defined on $\Sigma$ is cooperative, the set
$\mbox{Im }(\phi)$ must be unordered by Lemma~\ref{lemma unordered}, and
therefore $p^m\leq d_{m+1,p}$ by Lemma~\ref{lemma Sperner}. But by Lemma~\ref{dnpupper},
$d_{m+1,p} < p^m$, a contradiction.
\epr

Note that we are restricting our attention to the case where $p=q$, i.e. both
systems have the same level of discretization.   This is relevant for instance
in the special case of Boolean networks.  But if we allow $q\neq p$, then the
analogue of Theorem~\ref{theorem bounds}.3 may fail.

One important consequence of Theorem~\ref{theorem bounds} is that cooperative systems may have
exponentially long cycles, which can be considered a form of chaotic behavior in discrete systems.

\begin{Corollary}\label{chaoscorol}
Let $p > 1$ and let $c$ be an arbitrary real number with $1 < c < p$.  Then for sufficiently
large $n$, there exist $n$-dimensional cooperative $p$-discrete systems with periodic orbits of length
$> c^n$.
\end{Corollary}

\bpr
By Proposition~\ref{proposition embedding}, for each $m$ there exist $m$-dimensional cooperative $p$-discrete systems
with periodic orbits of length $d_{m,p}$.  Now the conclusion follows from Lemma~\ref{expolemma}.
\epr

While Lemmas~\ref{lemma dnp lower bound}-\ref{dnpupper} are sufficient for deriving our conclusions about
embeddings into cooperative systems, for completeness we will conclude this section with some sharper estimates of
$d_{n,p}$.
  The following result is an application of a local
central limit theorem.

\begin{Proposition} \label{proposition local limit}
For arbitrary $p>1$ and $\sigma^2=\frac{1}{12}(p-1)(p+1)$:
\[
\lim_{n\to \infty} d_{n,p}\  \left(\frac{p^{n}}{\sqrt{2\pi n
\sigma^2}}\right)^{-1}=1,
\]
\end{Proposition}

\bpr Let $Y_1,\ldots , Y_n$ be i.i.d. random variables, each of which can take
any of the values $0,\ldots, p-1$ with equal probability $1/p$.  Then
$E(Y_i)=(p-1)/2$, $\sigma^2(Y_i)=\frac{1}{12}(p-1)(p+1)$, for every $i$, and
$d_{n,p}=p^n \cdot P(\Sigma_{i=1}^n Y_i=\flo{ n (p-1) /2})$.

Define $X_i:=Y_i-(p-1)/2, \ S_n=\Sigma_{i=1}^n X_i$.  We will use the notation
of Section~2.5 in~\cite{Durrett:textbook}, in order to use Theorem~2.5.2 in
that textbook.  Set

\[
x=:\left(\flo{ \frac{n(p-1)}{2}} - \frac{n(p-1)}{2}\right)\frac{1}{\sqrt{n}}, \
\ \ p_n(x):=P(S_n/\sqrt{n}=x),
\]
and note that $p_n(x)=P(\Sigma_{i=1}^n Y_i - n(p-1)/2 = \sqrt{n}x) = d_{n,p}
p^{-n}$.

Let $\Phi(x):=(2\pi\sigma^2)^{-1/2}e^{-x^2/(2\sigma^2)}$. Then according to
Theorem~2.5.2, $\norma{\sqrt{n}p_n(x) - \Phi(x)}\to 0$ as $n\to \infty$.   But
evidently $\Phi(x)\to \Phi(0)$, since $x=x(n)\to 0$.  Therefore
$\sqrt{n}p_n(x)=\sqrt{n}d_{n,p}p^{-n} \to \Phi(0)$.  This immediately implies
the result. \epr

\begin{Corollary} \label{corollary dnp upper bound}
There exists constants $c_1, c_2>0$, independent of $n$ and $p$, such that for
 arbitrary $p>1$ there exists $n_0(p)$ such that for all $n\geq n_0(p)$:
\[
c_1 \frac{p^{n-1}}{\sqrt{n}}\leq d_{n,p} \leq  c_2 \frac{p^{n-1}}{\sqrt{n}}
\]
\end{Corollary}

\bpr Let $\sigma^2(p)=\frac{1}{12}(p-1)(p+1)$ be as in the above result. One
can compute $\sigma(p)/p=\frac{1}{\sqrt{12}}\sqrt{1-p^{-2}}$.  Therefore
$\frac{1}{4} p\leq \sigma(p) \leq \frac{1}{\sqrt{12}} p$, for all $p>1$.  Let
$\epsilon>0$ be an arbitrarily small number.  Then for all $n\geq
n_0=n_0(p)$
 Proposition~\ref{proposition local limit} implies:
\[
d_{n,p}\leq (1+\epsilon) \frac{ p^n}{\sqrt{2\pi n} \sigma(p)} \leq
\frac{1+\epsilon}{\frac{1}{4} \sqrt{2\pi}} \frac{p^{n-1}}{\sqrt{n}}.
\]
The second inequality is obtained analogously. \epr

\section{Smale extensions} \label{smaleapprox}

Assume $\Pi = \prod_{i=1}^n L_i$ and each $L_i$ has a smallest and a largest element.
Let $A \subseteq L$ and let $\gamma: A \rightarrow L$.  We say that $\gamma$ is \emph{cooperative} if for
all $x, y \in A$ the implication $x \leq y \rightarrow \gamma(x) \leq \gamma(y)$ holds. Clearly, if $A$ is unordered, then $\gamma$ is cooperative,
and since $L$ has a largest and a smallest
element, the construction used in the proof of
Lemma~\ref{lemma trivial embedding} implies that $\gamma$ can be extended to a cooperative
function on $\Pi$.  However, the construction used in this proof
is too crude to allow for such extensions if $A$ contains comparable elements.
Here we use a different construction to show that any cooperative partial function
on $\Pi$ can be extended to a cooperative function on $\Pi$.

\begin{Lemma}\label{smaleextension}
Let $\Pi = \prod_{i=1}^n L_i$, where $(L_i, <)$ is complete in the sense that every subset of $L_i$ has a supremum and an infimum in $L_i$.
Let $A \subseteq \Pi$,  and let
$\gamma: A \rightarrow \Pi$ be cooperative.  The there exists a
cooperative $g: \Pi \rightarrow \Pi$ such that $\gamma = g \rest A$.
\end{Lemma}

\bpr
Let $\Pi, A, \gamma$ be in the assumption.  First note that we may wlog assume that
for all $z \in \Pi$ there exists $x \in A$ such that $x \leq z$ or $z \leq x$.
If not, then extend $A$ to a set $A^*$ with this property and such that
$A^* \backslash A$ is unordered and each $x \in A^* \backslash A$ is incomparable
with each $z \in A$.  Extend $\gamma$ to $\gamma^*: A^* \rightarrow \Pi$ in an arbitrary
way, and note that $\gamma^*$ must still be cooperative.

Given $z\in \Pi$,  define
$U(z):=\{x\in A: \, x\geq z\}$
and $\Pi_U := \{z \in \Pi: \, U(z) \neq \emptyset\}$.
Note that $A \subseteq \Pi_U$.  Let $\Pi_L := \Pi \backslash \Pi_U$, and for all $z \in \Pi_L$ define
$L(z) :=  \{x\in \Pi_U: \, x\leq z\}$.  Note that our assumption on
 $A$ implies that $L(z) \neq \emptyset$ for all $z \in \Pi_L$.
Let $\gamma(U(z)) := \{\gamma(x): \, x \in U(z)\}$.

Now
define $g(z) :=
\inf \gamma(U(z))$ for $z \in \Pi_U$  and let
$g(L(z)) := \{g(x): \, x \in L(z)\}$ for $z \in L(z)$.  Finally, define $g(z) :=
\sup g(L(z))$ for $z \in \Pi_L$.

\begin{Claim} \label{prop Smale discrete}
The map $g$ defined above is cooperative and satisfies
$g \rest A=\gamma$.
\end{Claim}

\bpr
By completeness of $\leq$ on each $L_i$, infima and suprema under the cooperative order of nonempty subsets of $\Pi$ exist and are elements of $\Pi$.
Thus $g$ is well defined.

Suppose that $z\in A \subseteq \Pi_U$.  Then $\gamma(z) \leq \gamma(x)$ for every $x \in U(z)$ by cooperativity of $\gamma$, hence  $\gamma(z)  = \inf \gamma(U(z)) = g(z)$.
Thus $g \rest A = \gamma$.

To see that $g$ is
cooperative,  let $y,z\in \Pi$ be such that $y\leq z$.
If $y,z \in \Pi_U$, then $U(z)\subseteq U(y)$, and hence $g(y)=\inf \gamma(U(y))\leq \inf
\gamma(U(z))=g(z)$.
If $y,z \in \Pi_L$, then $L(y)\subseteq L(z)$, and hence
 then $g(y)=\sup g(L(y))\leq \sup
g(L(z)) =g(z)$. The only other possibility consistent with $y \leq z$ is
$z \in \Pi_L$ and $y \in \Pi_U$.
In this case $y \in L(z)$, and hence $g(y) \leq \sup g(L(z)) = g(z)$. \epr

\epr

For reasons that will become clear shortly, we will refer to the function $g$ constructed in the proof of
Lemma~\ref{smaleextension} as the \emph{Smale extension of $\gamma$.}

Now suppose $\Pi$ is is either $\{0, \ldots , p-1\}^n$ or $[0,1]^n$ with the natural order, and let $A$
be a hyperplane of the form $A = \{x \in \Pi: \, S(x) = r\}$. If $A$ is nonempty (which will happen for
suitable values of $r$), then $A$ is a maximal incomparable subset of $\Pi$.

Note that if $A$ is a hyperplane as above, then the definition of the Smale extension $g$ of~$\gamma$ can be written as

\begin{equation}\label{altsmaledef}
g(z)= \left\{
\begin{array}{ll}
\inf \gamma(U(z)), & \mbox{ for } z \in \Pi_U,  \\
\sup \gamma(L(z)), &  \mbox{ for } z \in \Pi_L.
\end{array} \right.
\end{equation}

 For $x \in \Pi$, let $\| x \| = \max \{|x_1|, \ldots , |x_n|\}$ be the sup-norm in ${\R}^n$.

\begin{Lemma}\label{epsdeltalemma}
Suppose $\Pi$ is  either $\{0, \ldots , p-1\}^n$ or $[0,1]^n$ with the natural order, and
$A= \{x \in \Pi: \, S(x) = r\}$ is a nonempty hyperplane.
 Let $\gamma, \gamma_1 : A \rightarrow \Pi$ be cooperative, and let $\varepsilon, \delta > 0$ be such that
\begin{equation}\label{epsdeltaeqn}
\forall x, y \in A \quad \|x - y\| < (2n+1) \delta \Rightarrow \|\gamma (x) - \gamma_1 (y)\| < \frac{\varepsilon}{3}.
\end{equation}

Let $g, g_1$ be the Smale extensions of $\gamma, \gamma_1$.  Then

\begin{equation}\label{epsdeltaeqn1}
\forall x, y \in \Pi \quad \|x - y\| < \delta \Rightarrow \|g (x) - g_1 (y)\| < \varepsilon.
\end{equation}

\end{Lemma}

\bpr
Let $A, \gamma, \varepsilon, \delta$ be as in the assumption.  First note that

\begin{equation}\label{closeb}
\forall x, y \in \Pi_U \forall a \in U(x)\exists b \in U(y) \  \|a - b\| \leq (n-1)\|x - y\|.
\end{equation}

To see this, let $x, y \in \Pi_U$ and $a \in A$ with $x \leq a$.  Let $b \in U(y)$ be such that $\sum |b_i - a_i|$ is minimal.
Such $b$ exists by compactness of $A$.
 Note that we must have
$b_i \leq \max \{a_i, y_i\}$ for all $i \in \{1, \ldots , n\}$: If not, since $S(a) = S(b)$, there would be some $j$ with
$y_j \leq b_j < a_j$.  Letting $\beta = \min \{b_i - \max \{a_i, y_i\}, a_j - b_j\}$ and $b^*_i = b_i - \beta$, $b^*_j = b_j + \beta$,
and $b^*_k = b_k$ for all $k \neq i, j$,
we would have $y \leq b^* \in A$ and $\sum |b_i^* - a_i| < \sum |b_i - a_i|$, contradicting the choice of $b$.  Thus
$|b_i - a_i| \leq \|x - y\|$ for all $i$ with $b_i > a_i \geq x_i$.  Now consider $i$ with $y_i \leq b_i < a_i$.  In this situation we must
have $a_i - b_i \leq \sum \max \{0, b_j - a_j\} \leq (n-1) \|x - y\|$.  The inequality $\|a - b \| \leq (n-1)\|x - y\|$ follows.

Furthermore, note that

\begin{equation}\label{closeb1}
\forall x \in \Pi_U \forall y \in \Pi_L \forall a \in U(x) \forall b \in L(y) \,
 \|a - b\| \leq (2n+1)\|x - y\|.
\end{equation}

To see this, let $x \in \Pi_U, y \in \Pi_L,  a \in U(x), b \in L(y)$.  Fix $i \in \{1, \ldots , n\}$.  Then $a_i - x_i \leq S(a) - S(x) =
r - S(x) \leq S(y) - S(x) \leq n\|y - x\|$.  Similarly,  $y_i - b_i \leq S(y) - S(b) =
S(y) - r \leq S(y) - S(x)$.  Now it follows from the triangle inequality that $|a_i - b_i| \leq 2(S(y)-S(x)) + |y_i - x_i| \leq
(2n+1)\|y - x\|$.

Now let $x, y$ be such that $\|x - y\| < \delta$.

First  assume that $x, y \in \Pi_U$.
Fix $i \in n$, and let $a \in U(x)$ be such that  $|(g(x))_i - (\gamma(a))_i| < \varepsilon/3$.  Such $a$ exists by~(\ref{altsmaledef}).  Choose $b\in U(y)$ as in~(\ref{closeb}).
It follows from~(\ref{epsdeltaeqn})  that $\|\gamma (a) - \gamma_1 (b)\| < \varepsilon/3$.  In particular, $|(\gamma(a))_i - \gamma_1(b)_i| < \varepsilon/3$. Since $y \leq b$, definition~(\ref{altsmaledef}) implies that
$(g_1(y)_i) \leq (\gamma_1(b))_i$, and the inequality
 $(g_1(y))_i < (g(x))_i + 2 \varepsilon/3$ follows.  By symmetry of the assumption, we also will have $(g(x))_i < (g_1(y))_i + 2 \varepsilon/3$ in this case.

By the alternative definition~(\ref{altsmaledef}) of the Smale embedding, the argument in the case when  $x, y \in \Pi_L$ is dual.

Now assume $x \in \Pi_U$ and $y \in \Pi_L$.  Fix $i \in n$, and let $a \in U(a)$ and $b \in L(b)$ be such that
$|(g(x))_i - (\gamma(a))_i| < \varepsilon/3$ and $|(g_1(y))_i - (\gamma_1(b))_i| < \varepsilon/3$.  By~(\ref{closeb1}),
$\|a - b\| \leq
(2n+1)\|y - x\| <(2n+1)\delta$, and~(\ref{epsdeltaeqn}) implies that $|(\gamma(a))_i - \gamma_1(b)_i| < \varepsilon$. Now~(\ref{epsdeltaeqn1})
follows from the triangle inequality.

The argument in the case when  $x \in \Pi_L$ and $y \in \Pi_U$ is symmetric.
\epr

By letting $\gamma = \gamma_1$ in Lemma~\ref{epsdeltalemma}, we immediately get the following:

\begin{Corollary}\label{Lipshitzcorol}
Suppose $\Pi = [0,1]^n$ with the natural cooperative order, and
$A= \{x \in \Pi: \, S(x) = r\}$ for some $0 \leq r \leq 1$.
 Let $\gamma: A \rightarrow \Pi$ be cooperative, and let $g: \Pi \rightarrow \Pi$ be the Smale extension of $\gamma$.

\noindent
(i) If $\gamma$ is continuous, so is $g$.

\noindent
(ii) If $\gamma$ is Lipshitz-continuous with Lipshitz constant $\ell$, then $g$ is Lipshitz continuous with Lipshitz
constant $\leq (6n+3)\ell$.
\end{Corollary}

Now consider any discrete-time  dynamical
system  $([0,1]^n, f)$, let $A = \{x \in [0,1]^{n+1}: S(x) = (n+1)/2\}$, and let $\phi: [0,1]^n \rightarrow A$ be a Lipshitz-continuous homeomorphism.
Let $\gamma: A \rightarrow A$ be such that $\gamma(\phi(x)) = \phi(f(x))$ for all $x \in [0,1]$.
If $f$ is (Lipshitz)-continuous, then so is $\gamma$, and Corollary~\ref{Lipshitzcorol} implies that $\phi$ is a
(Lipshitz)-continuous embedding of $([0,1]^n, f)$ into a
discrete-time  dynamical
system  $([0,1]^{n+1}, g)$ for which $g$ is (Lipshitz)-continuous.
This is analogous to Smale's famous embedding theorem for $C^1$-systems \cite{Smale:1976}
and is our motivation for calling the function~$g$ of Lemma~\ref{smaleextension} the \emph{Smale extension} of $\gamma$.

If $([0,1]^n, f)$ and $(\{0, \ldots p-1\}^n, f_1)$ are two discrete-time systems and $\varepsilon > 0$, then
we will say that \emph{$f_1$ is an $\varepsilon$-approximation of $f$} if $\|\frac{1}{p-1}f_1(\lfloor(p-1)x\rfloor) - f(x)\| < \varepsilon$ for all $x \in [0,1]^n$.

Let $A$ be as in the previous paragraph,  let $D = \{y \in \{0, \ldots , p-1\}^{n+1}: \,
S(y) = \lfloor(n+1)(p-1)/2\rfloor\}$, and define $D^* := \{a \in A: \, (p-1)a \in D\}$. It is clear that if $f$ is continuous,
$\delta > 0$ is given, $\beta > 0$ is sufficiently small relative to $\delta$,
$p$ is odd and sufficiently large,
and if $\phi, \gamma$ are as in the previous paragraph, then there exist:
\begin{itemize}
\item [-] a $\beta$-approximation $(\{0, \ldots , p-1\}^n, f_1)$
of $([0,1]^n, f)$,
\item [-] and a function $\gamma_1: A \rightarrow D^*$ such that $\|\gamma(y) - \gamma_1(y)\| < \delta$ for all
$y \in A$,
\item [-] a function $\phi_1: \{0, \ldots p-1\}^n \rightarrow D^*$
such that
$\|\phi(x/(p-1)) - \phi_1(x)\| < \delta$ and $\gamma_1(\phi_1(x)) = \phi_1(f(x))\|$ for all $x \in \{0, \ldots , p-1\}^n$.
\end{itemize}

Now let $\varepsilon > 0$ be given. By Lemma~\ref{epsdeltalemma}, if we choose the above objects for $\delta$ sufficiently small relative to $\varepsilon$ and
 if  $g$ is the Smale extension of $\gamma$, while $g_1$ is the Smale extension of
$\gamma_1$, then $\|g(y) - g_1(y)\| <  \varepsilon$ for all $y \in [0,1]^{n+1}$.
Let $g_1^*(x) := (p-1)g_1(x/(p-1))$ for all $x \in \{0, \ldots , p-1\}^{n+1}$.
From the definition of the Smale extension it follows that $g_1^*$ maps $x \in \{0, \ldots , p-1\}^{n+1}$ into itself, and
the inequality $\|g(y) - g_1(y)\| <  \varepsilon$ implies that $g^*_1$ is an $\varepsilon$-approximation of $g$.
Moreover, we will have $g^*_1(\phi_1(x)) = \phi_1(f_1(x))$ for all
$x \in \{1, \ldots ,1\}^{n}$.  However, we cannot necessarily assume that $\phi_1$ is a cooperative embedding
of $(\{0, \ldots p-1\}^n, f_1)$ into $(\{0, \ldots , p-1\}^{n+1}, g^*_1)$ , since the results of Section~\ref{section bounds}
indicate that the function $\phi_1$ may not be injective.

\section{Almost cooperative systems}\label{almost coop section}

Cooperative systems are so named because increasing the value of one variable
tends to increase the values of other variables in the system.
For instance, in the continuous case a condition equivalent to the
cooperativity of the system (\ref{ode}) is $\partial f_i / \partial x_j (x)\geq
0$ for $i\not=j$ \cite{Smith:monotone}.
 It has been conjectured that a system might have amenable properties if it is `almost cooperative,' i.e. if the latter condition is satisfied with the exception of a single pair $i\not=j$ (see the concept of \emph{consistency deficit} in \cite{Sontag:SSB2007}).

We can define a discrete counterpart of this notion as follows.
Let $x \in \Pi = \prod_{i=1}^n\{0, \ldots , p_i-1\}^n$,
and let $i \in \{1, \ldots , n\}$.  Define $x^{i+} \in \Pi$ by
letting $(x^{i+})_i = \min \{x_i + 1, p_i-1\}$ and $(x^{i+})_j = x_j$ for $j \neq i$.  Similarly,
define $x^{i-} \in \Pi$ by
letting $(x^{i-})_i = \max \{x_i - 1, 0\}$ and $(x^{i-})_j = x_j$ for $j \neq i$.
It is easy to see that cooperativity of a system $(\Pi, g)$ is equivalent to the condition
that

\begin{equation}\label{coopi+def}
\forall x \in \Pi \quad (g(x^{i-}))_j \leq (g(x))_j \leq (g(x^{i+}))_j
\end{equation}
for all $i, j \in \{1, \ldots , n\}$.  Let us call $(\Pi, g)$ \emph{almost cooperative}
if condition~(\ref{coopi+def}) holds with the exception of exactly one
pair $<i^*, j^*>$ with $i^* \neq j^*$ for which an order-reversal takes place:
\begin{equation}\label{compi+def}
\forall x \in \Pi \quad (g(x^{i^*-}))_{j^*} \geq (g(x))_{j^*} \geq (g(x^{i^*+}))_{j^*}.
\end{equation}

 One might expect that almost cooperative
$p$-discrete systems are similar to cooperative systems. In particular,
one might expect that  $m$-dimensional almost cooperative Boolean systems
 can always be
embedded into  cooperative Boolean systems of dimension $m+1$.  However, this is
not the case.

Consider the following simple
example with $n=2$.  In this case $|\Pi| = 4$, and $d_{3,2}=\binom{3}{2}=3$.
Define $g(x_1, x_2) := (1-x_2 , x_1)$, so that $g(0,0) = (1,0)$, $g(1, 0) =
(1,1)$, $g(1,1) = (0,1)$, and $g(0,1) = (0,0)$.  The system  $(\Pi, g)$
consists of a single orbit of length $4$.  By
Proposition~\ref{proposition embedding}, this system cannot be embedded into
any cooperative Boolean system of dimension $3$.   Moreover, note that $x_1$
promotes the increase of the variable $x_2$, while $x_2$ inhibits the variable
$x_1$.  Thus condition~(\ref{coopi+def}) holds with the exception of the pair
$<i,j>\ =\ <2,1>$ and hence the
system is almost cooperative.

\section{Strong Cooperativity for Finite Discrete Systems}  \label{section alternative}

Throughout the remaining two sections we will assume that the state space $\Pi$ of our dynamical system is
finite.  Wlog this means that $\Pi = \prod_{i=1}^n \{0, \ldots , p_i -1\}$,
where the $p_i$'s are integers such that $p_1 \geq \dots \geq p_i \geq \dots \geq p_n > 1$.
Of course, $p$-discrete systems are exactly those among the above systems for which $p_i = p >1$
for all~$i$.

Our first task is to come up with a suitable counterpart of strong
cooperativity for such systems.  Let us write $x < y$ if $x \leq y$ in the cooperative order but
$x \neq y$, and let us write $x \ll y$ if $x_i < y_i$ for all $i = 1, \ldots ,
n$. Recall that a continuous system is strongly cooperative if for every two
initial conditions $x(0) < y(0)$ we have $x(t) \ll y(t)$ for all $t
> 0$.
Perhaps the most straightforward adaptation of this definition to finite discrete
systems would be the following:

\begin{equation}\label{sm0def}
\forall x(0) < y(0) \, \exists t_0 > 0 \, \forall t \geq t_0\quad x(t) \ll
y(t).
\end{equation}
This property is commonly known as eventual strong cooperativity
\cite{Smith:monotone}.  Unfortunately,  no finite discrete system of dimension $n
> 1$ satisfies~(\ref{sm0def}).  To see this, note that if $n > 1$, then there
exist $x^0(0) < x^1(0) < \dots < x^{p_1}(0)$.  On the other hand, $x(t) \ll y(t)$
implies $n + S(x(t)) \leq S(y(t))$.  Now~(\ref{sm0def}) would imply
$S(x^{p_1}(t)) \geq np_1$ for sufficiently large $t$. But this is a contradiction,
since $S( x) \leq (p_1 -1)n$ for all $x \in \Pi$.

So let us define a weaker discrete version of strong cooperativity.
Consider the following properties of a finite cooperative system.

\begin{equation}\label{sm1def}
\forall x(0) < y(0) \quad  S (y(0)- x(0)) \leq  S(y(1)- x(1)).
\end{equation}

\begin{equation}\label{sm2def}
\forall x(0) < y(0) \quad   x(1) < y(1).
\end{equation}

\begin{equation}\label{smsumdef}
\forall x(0) \quad   S (x(1)) = S(x(0)).
\end{equation}

\begin{Lemma}\label{sm12lemma}
For any finite cooperative discrete system conditions~(\ref{sm1def}),~(\ref{sm2def})
and~(\ref{smsumdef}) are equivalent.
\end{Lemma}

\bpr Clearly, condition~(\ref{sm1def}) implies condition~(\ref{sm2def}) and
condition~(\ref{smsumdef}) implies condition~(\ref{sm1def}) for cooperative systems.

Now assume an $n$-dimensional finite system $(\Pi, g)$ satisfies condition
(\ref{sm2def}).  Let $x(0) \in \Pi$, let $S := \sum_{i=1}^n (p_i - 1)$, and consider initial states $x^0(0) <
x^1(0) < \dots < x^{S}(0)$ such that $x(0) = x^{S(x)}(0)$.
 By (\ref{sm2def}),  $x^0(1) < x^1(1) < \dots <
x^{S}(1)$.
It follows that $S(x^i(1))=S(x^i(0))$ for all $i$, and in particular
(\ref{smsumdef}) holds.
 \epr

We will say that a finite cooperative system $(\Pi, g)$ is \emph{strongly
cooperative} if it satisfies conditions~(\ref{sm1def}),~(\ref{sm2def})
and~(\ref{smsumdef}).

For example,  let $\pi: \{1, \ldots , n\} \rightarrow \{1, \ldots n\}$ be a
permutation, and define $g_\pi$ by $(g_\pi(x))_{\pi(i)} = x_i$ for all
$x \in \Pi$   and
$i \in \{1, \ldots , n\}$.   Note that if $\Pi = \{0, \ldots , p-1\}^n$, then $g_\pi$ maps  $\Pi$ into $\Pi$ for all permutations
 $\pi$ of $\{1, \ldots , n\}$, but in general this will be the case  only for some but not for all
 permutations. If $g_\pi$ does map $\Pi$ into $\Pi$, then $(\Pi, g_\pi)$ is a strongly cooperative
 system.

The \emph{order} of a permutation $\pi$ is the smallest
integer $r > 0$ such that $\pi^r$ is the identity. Let $R(n)$ be the maximum order of a permutation $\pi$ of $\{1, \ldots , n\}$.  It can be shown that
$R(n)= e^{ \sqrt{n \ln n} (1+o(1)) }$  as $n\to \infty$ \cite{Landau}.
%gae4 that $R(n) < c2^{\sqrt{n}}$ for some constant $c > 0$.
%\textbf{Get $c$ from literature.}
In particular, note that $R(n)$ grows \emph{subexponentially in $n$,}
that is, for every $b > 1$ and sufficiently large $n$ we will have $R(n) < b^n$.

\begin{Theorem}\label{smcyclelem}
Suppose $( \prod_{i=1}^n \{1, \ldots p_i - 1\}^n, g)$ is an $n$-dimensional strongly cooperative finite discrete system,
and let $N = \sum_{i=1}^n (p_i -1)$.
Then each periodic orbit in $( \Pi, g)$ has length at most  $R(N)$.
\end{Theorem}

\bpr Note that for any permutation $\pi$, the length of
any periodic orbit of $(\Pi, g_\pi)$ cannot exceed the order of $\pi$.
However, not all strongly cooperative finite systems are of the form
 $(\Pi, g_\pi)$ for some permutation $\pi$.
For example, if $n = 2$ and $g(\emptyset) = \emptyset$, $g (\{1\}) = g(\{2\}) =
\{1\}$ and $g(\{1,2\}) = \{1,2\}$, then $g$ is a strongly cooperative Boolean system, but not of
the form $g_\pi$ for any permutation~$\pi$.

Fortunately, Lemma~\ref{piclaim} below
suffices for the proof of our theorem in the Boolean case when $p_i = 2$ for all $i$ and hence $N = n$.  A state in an attractor of
a dynamical system, i.e., a state that is not transient, will be called
\emph{persistent.}

\begin{Lemma}\label{piclaim}
Let $(\Pi,g)$ be a strongly cooperative $n$-dimensional Boolean system.  Then
there exists a permutation~$\pi$ of $\{1, \ldots , n\}$ such that $g(x) =
g_\pi(x)$  for each persistent state $x$ of $(\Pi, g)$.
\end{Lemma}

\bpr
We will prove the lemma by induction over $n$.  Note that it is trivially
true for $n = 1$.  Now fix $n > 1$, assume the lemma is true for all $k < n$
and let $(\Pi , g)$ be as in the assumption. We will identify
elements $x$ of $\Pi$ with subsets of the set $\{1, \ldots , n\}$ and write
$|x|$ instead of $S(x)$.  Note that $g$ maps one-element
subsets of $\{1, \ldots ,n\}$ to one-element subsets.  More precisely, there
exists a function $\sigma : \{1, \ldots , n\} \rightarrow \{1, \ldots , n\}$
such that $g(\{i\}) = \{ \sigma(i)\}$ for all $i$.  In general, $\sigma$ does
not need to be a bijection.  However, if $I$ is the set of all $i$ such that
$\{i\}$ is a persistent state of our system, then $I \neq \emptyset$ and
$\sigma \upharpoonright I$ is a permutation of $I$. Now strong cooperativity of
$g$ implies that $g(x) = g_{\sigma \rest I}(x)$ for all $x \subseteq I$. Thus
if $I = \{1, \ldots , n\}$, we are done.  If not, then define for $y \in \Pi$
such that $I \cap y = \emptyset$:

$$f(y) = g(y \cup I) \backslash I.$$

Since $g(I) = g_{\sigma \rest I}(I) = I$, the function $f$ is strongly
cooperative on the set of all subsets of $J := \{1, \ldots , n\} \backslash I$.
By the inductive assumption, there exists a permutation $\varrho$ of $J$ such
that $f(y) = g_\varrho(y)$ for all persistent states in the system defined by
$f$.  Note that $\pi := (\sigma \rest I) \cup \varrho$ is a permutation of
$\{1, \ldots , n\}$.

Now consider any $x = x(0) \in \Pi$.  By strong cooperativity we have
$$|x| = |g(x)| = |g(x) \cap I| + |g(x) \cap J|.$$

On the other hand, $|x \cap I| = |g(x \cap I)| \leq |g(x) \cap I|$ because $g(x
\cap I) \subseteq g(I) = I$.  It follows that $|x(t) \cap J|$ is nonincreasing
along the trajectory of $x(0)$.  In particular, for every persistent state $x$
we must have $g(x \cap I) = g(x) \cap I$ and hence $|x \cap J| = |g(x) \cap J|$.

It must also be the case that $g(x) \cap J \subseteq f(x \cap J)$.  Since $|f(x \cap J)| =
|x \cap J|$ by strong cooperativity of $f$, we must have  $g(x) \cap J = f(x
\cap J)$ for every persistent state $x$ of $(\Pi, g)$.  It follows that if $x$
is a persistent state of $(\Pi, g)$, then $x \cap J$ is a persistent state of
$(\Sigma, f)$.  Thus $g(x) = g_{\sigma \rest I} (x \cap I) \cup g_\varrho(x
\cap J) = g_\pi(x)$.
\epr

Now consider the general case where $p_i \geq 2$ for all $i$.  Unfortunately, we cannot hope to
prove the exact analogue of Lemma~\ref{piclaim}.
To see this, consider the system
$(\{0,1,2\}^2, g)$, where $g(0,0) = [0, 0]$, $g(0,1) = [1,0]$, $g(1, 0) =
[0,1]$, $g(1,1) = g(2,0) = (0,2) = [1,1]$, $g(1,2) = [1,2]$, $g[2,1] = [2,1]$, $g(2,2) = [2,2]$.
  This system is clearly strongly
cooperative.  If $\pi$ were a permutation as in Lemma~\ref{piclaim}, then we would
need $\pi(0) = 1$ and $\pi(1) = 0$ because both $[0,1]$ and $[1,0]$ are persistent states.
On the other hand, $[2,1]$ and $[1,2]$ are persistent steady states, so this would force
 $\pi$ to be the identity.

Now let $\Pi = \prod_{i=1}^n \{1, \ldots p_i - 1\}^n$ and let $(\Pi, g)$, $N$ be as in the assumption of
Theorem~\ref{smcyclelem}.  Let $\Sigma :=
\{0, 1\}^N$.
For each $i \in \{1, \ldots , n\}$ let $J_i$ be the set of integers $j$ such that
$\sum_{k=1}^{i-1} (p_k-1) < j \leq \sum_{k=1}^{i-1} (p_k-1)$. For $1 \leq \ell \leq |J_i| = p_i - 1$ let $j(i, \ell)$ be the
$\ell$-th element of $J_i$.  Define a map $\psi: \Pi \rightarrow \Sigma$ so that
for  $i \in \{1, \ldots , n\}$ and $\ell \in \{1, \ldots ,  p_i-1\}$ we have
$\psi(x)_{j(i, \ell)} = 1$ iff $x_i \geq \ell$.  Clearly, $\psi$ is an
injection. For $y \in \Sigma$ define $z(y)$ by $z(y)_{j(i, \ell)} = 1$ iff
$\ell  \leq |\{\ell': \ y_{j(i, \ell')} = 1\}|$. Note that $z(y)$ is always in the
range of $\psi$, and $z(y) = y$ whenever $y$ is already in the range of $\psi$.
Moreover, the function $y \mapsto z(y)$ is strongly cooperative. Now define $f:
\Sigma \rightarrow \Sigma$ so that $f(y) = f(z(y))$ for all $y$ and $f(\psi(x)) =
\psi(g(x))$ for all $x \in \Pi$.

Then $\psi$ is an embedding of $(\Pi, g)$ into $(\Sigma, f)$.  Moreover, if
$(\Pi, g)$ is strongly cooperative, then so is $(\Sigma, f)$.
 Since the lemma is true for $p = 2$, each periodic orbit in  $(\Sigma, f)$ has length at most $R(N)$,
 and since $\psi$ is an embedding, the same must be true for $(\Pi, g)$.
\epr

\begin{Question}\label{orderquestion}
Suppose $(\Pi, g)$ is an arbitrary $n$-dimensional strongly cooperative finite discrete system.
Can the system have a periodic orbit of length greater than  $R(n)$?  What if we assume in
addition that $(\Pi, g)$ is $p$-discrete for some $p > 2$?
\end{Question}

Our results can perhaps be considered analogues of the result in
\cite{Terescak:1992} for discrete-time continuous-space strongly cooperative
systems. Our Theorem~\ref{smcyclelem} gives a nontrival, subexponential bound on
the lengths of periodic orbits of strongly cooperative finite discrete systems. Moreover,
 strong cooperativity implies that ordered orbits in finite discrete systems are fairly
robust, as shown in the next result.

\begin{Lemma}  Consider a strongly cooperative finite discrete system (\ref{discrete}), and let
$x(0)$ and $y(0)$ be two arbitrary initial conditions (i.e.\ not necessarily ordered).  Then $S(|y(t) -
x(t)|) \leq S(|y(0) - x(0)|) $ for all $t > 0$.
\end{Lemma}

\bpr
Suppose first that $S(|y(0) - x(0)|)=1$.  Then necessarily the two initial conditions are ordered; suppose without loss of generality $x(0)<y(0)$.   By condition~(\ref{smsumdef}) and cooperativity, $S(|y(0) - x(0)|)=S(y(0) - x(0))=S(y(0))-S(x(0))=S(y(1))-S(x(1))= S(|y(1) - x(1)|)$.

If $S(|y(0) - x(0)|)=k>1$, then there exists a sequence of states $x=x^0,x^1,\ldots,x^k=y$, such that $S(|x^{j+1}-x^j|)=1$ for every $j$.  Then $S(|y(1) - x(1)|) \leq S(|x^k(1) - x^{k-1}(1)|)+\ldots + S(|x^1(1) - x^0(1)|)= S(|x^k(0) - x^{k-1}(0)|)+\ldots + S(|x^1(0) - x^0(0)|)=k=S(|y(0) - x(0)|)$.
\epr

In other words, small perturbations of initial
conditions don't amplify along the trajectory, which implies an analogue of
Lyapunov stability for all attractors.

\section{Cooperative Irreducible Systems and Long Periodic Orbits}  \label{section strong coop}

In this section we will explore several possible discrete counterparts of the notion of irreducible
cooperative $C^1$-systems and will show how these conditions relate to strong cooperativity and what bounds they
impose on the lengths of periodic orbits.

Recall that a digraph (directed graph) $G = (V, A)$ is \emph{strongly connected} if every node
$w$ in $V$ can be reached via a directed path from every node $v \in V$.

Now let us define discrete analogues of irreducible cooperative systems by
associating directed graphs $G = (\{1, \ldots , n\}, A)$ with a cooperative
system $(\Pi, g)$. Recall that in the definition of irreducible cooperative
$C^1$-systems, an arc $<i,j>$ was included in the arc set $A$ iff
$D f(x)_{ij}>0$ on $\R^n$, where $D f(x)$ is the Jacobian of $f(x)$, and the system was called irreducible if
the resulting directed graph $G$ on $\R^n$ was
strongly connected.  Alternatively, a digraph $G_x$ can be defined locally for every
$x\in \R^n$ by letting $<i,j>$ be an arc in $G_x$ if and only if $D f(x)_{ij}>0$.
  A cooperative $C^1$-system in which $G_x$ is strongly connected for every
$x \in \R^n$ is still strongly cooperative;  see for instance Corollary~3.11 in \cite{Smith:review}.

Recall the definitions of $x^{i-}$ and $x^{i+}$ from Section~\ref{almost coop section}.
For an $n$-dimensional  finite discrete system $(\Pi, g)$ and $x \in \Pi$, let us define a directed graph
$G^*_x = (\{1, \ldots , n\}, A^*_x)$ by including an arc $<i, j> \ \in A^*_x$ iff
$g(x)_j < g(x^{i+})_j$ or  $g(x^{i-})_j < g(x)_j$.
Moreover, let us define a directed graph
$G_x = (\{1, \ldots , n\}, A_x)$ by including an arc $<i, j> \ \in A_x$ iff
$<i,j> \ \in A^*_x$ and if $0 < x_i < p_i -1$, then $g(x^{i-})_j < g(x)_j < g(x^{i+})_j$.

Let us call the system $(\Pi, g)$ \emph{strongly irreducible} if $(\{1, \ldots n\}, \bigcap_{x \in \Pi} A_x)$ is
strongly connected, \emph{strongly semi-irreducible} if $(\{1, \ldots n\}, \bigcap_{x \in \Pi} A^*_x)$ is
strongly connected, \emph{irreducible} if $G_x$ is
strongly connected for all $x \in \Pi$, and \emph{weakly irreducible} if $(\{1, \ldots n\}, \bigcup_{x \in \Pi} A^*_x)$ is
strongly connected.

Note that for Boolean systems, $A^*_x = A_x$ for all states $x$; hence the notions of strong irreducibility and
strong semi-irreducibility coincide for Boolean systems.  While both strong irreducibility and strong semi-irreducibility
are plausible counterparts of irreducibility in cooperative $C^1$-systems, we will see that these two notions have dramatically
different implications for the dynamics of non-Boolean finite discrete systems.

In analogy to continuous systems, one would expect that irreducibility of cooperative discrete systems would imply
strong cooperativity and would put nontrivial bounds
on the length of periodic orbits, but weak irreducibility would not, since this condition can be guaranteed in local neighborhoods
of $x$'s that are far away from the attractor.  Therefore we will also consider the following properties.
Let $(\Pi, g)$ be an $n$-dimensional finite discrete system, and let $X$ be an attractor.  We say that
the system $(\Pi, g)$ is \emph{strongly irreducible along $X$} if $(\{1, \ldots n\}, \bigcap_{x \in X} A_x)$ is
strongly connected, \emph{irreducible along $X$} if $G_x$ is
strongly connected for all $x \in X$, and \emph{weakly irreducible along $X$} if $(\{1, \ldots n\}, \bigcup_{x \in X} A^*_x)$ is
strongly connected.  Note that (strongly) irreducible systems are (strongly) irreducible along every attractor.  In contrast,
a system that is weakly irreducible along at least one attractor is already weakly irreducible.

Let $\pi$ be a permutation of $\{1, \ldots , n\}$.  Recall the definition of the function $g_\pi$
 from the previous section.  Note that a finite discrete system $(\Pi, g_\pi)$ is strongly irreducible iff
 $(\Pi, g_\pi)$ is weakly irreducible along some attractor iff
 the permutation $\pi$ is cyclic.

\begin{Theorem}\label{itosc}
Suppose $(\Pi, g)$ is a cooperative irreducible finite discrete system.  Then $(\Pi, g)$ is strongly cooperative and
strongly irreducible.  Moreover, there exists a cyclic permutation $\pi$ of $\{1, \ldots , n\}$ such that
$g = g_\pi$.
\end{Theorem}

\bpr
Let $x, y \in \Pi$ be such that $x < y$.  Pick  $i \in \{1, \ldots , n\}$ such that $x < x^{i+} \leq y$.  Then there exists some
$j$ with $<i, j> \ \in A_x$; otherwise $G_x$ could not be strongly connected.  Thus $g(x)_j < g(x^{i+})_j \leq g(y)$ by cooperativity,
and condition~(\ref{sm2def}) follows.  Thus $(\Pi, g)$ is strongly cooperative.

Now let us consider $A_{\vec{0}}$.  Let us write $\{i\}$ for the $x \in \Pi$ with $x_i = 1 = S(x)$ and call such $x$ a \emph{singleton.}
Note that $<i, j> \ \in A_{\vec{0}}$ iff $g(\{i\})_j > 0$.
By strong cooperativity, $g(\{i\})$ is a singleton, and it follows that the outdegree of each $i$ in $G_{\vec{0}}$ is at most one.
Strong connectedness of $G_{\vec{0}}$ now implies that the in- and outdegrees in $G_{\vec{0}}$ of all nodes are exactly one.
Let $\pi: \{1, \ldots , n\} \rightarrow \{ 1, \ldots , n\}$ be defined  by $\pi(i) = j$ iff $<i, j>\ \in A_{\vec{0}}$.  Then $\pi$ is
a permutation.  Moreover, if $\pi$ could be decomposed into nonempty pairwise disjoint cycles, then $G_{\vec{0}}$ would not be strongly
connected.  Thus $\pi$ must be cyclic.

It remains to show that $g = g_\pi$.  We will show this by induction over $S(x)$.   If $S(x) = 0$, then $g(x) = x$ by strong cooperativity,
hence $g(x) = g_\pi(x)$.  By the way we defined $\pi$ we also have
$g(x) = g_\pi(x)$ whenever $S(x) = 1$.

Now let us assume $g(x) = g_\pi(x)$ for all $x$ with $S( x) = k$, and let $y$ be such that $S( y) = k+1$.  Then $y = x^{i+}$ for some $i$ and
$x$ with $S(x) = k$.  By the inductive assumption, $g(x) = g_\pi(x)$.  If $x_i = 0$, then $(g_\pi(x))_{\pi(i)} = 0$ but we must have both
$g(x) < g(y)$ and $g(\{i\}) = \{j\} \leq g(y)$, so $g_\pi(y) \leq g(y)$, and strong cooperativity implies $g(y) = g_\pi(y)$.
If $x_i > 0$, then the definition of $A_x$ implies that there must be $j$ with $g(x^{i-})_j < g(x)_j < g(x^{i+})_j$.  But by inductive assumption,
the only $j$ with $g(x^{i-})_j < g(x)_j$ is $\pi(i)$, so we must also have $g(x)_j < g(x^{i+})_j$.  It again follows that
$g_\pi(y) \leq g(y)$, and hence  $g(y) = g_\pi(y)$ by strong cooperativity.
\epr

\begin{Corollary}\label{shrtcyccor}
Periodic orbits in cooperative irreducible $n$-dimensional finite discrete systems
can have length at most $n$.
\end{Corollary}

\bpr
The maximal order of a cyclic permutation on $\{1, \ldots , n\}$ is $n$.
\epr

Corollary~\ref{shrtcyccor} gives a stronger bound than Theorem~\ref{smcyclelem} does for strongly cooperative
$p$-discrete systems.  For Boolean systems we can prove the same bound under weaker assumptions.

\begin{Lemma}\label{wcscshort}
Suppose $(\Pi, g)$ is a strongly cooperative $n$-dimensional Boolean system, and let $X$ be an attractor such that $(\Pi, g)$ is weakly irreducible
along $X$.  Then $|X| \leq n$.
\end{Lemma}

\bpr
Let  $(\Pi, g)$ and $X$ be as in the assumptions.  We will identify $x \in \Pi$ with subsets of $\{1, \ldots , n\}$.
By Lemma~\ref{piclaim}, there exists a permutation $\pi$ of
$\{1, \ldots , n\}$, where $n$ is the dimension of the system, such that $g(x) = g_\pi(x)$ for all $x \in X$.
It suffices to show that $\pi$ is cyclic.  Let $I \subset \{1, \ldots , n\}$ be the elements of a cycle of
$\pi$ such that $x \cap I \neq \emptyset$ for all $x \in X$ or $X = \{\emptyset\}$.
Now suppose $ I \neq \{1, \ldots , n\}$, and let $J = \{1, \ldots , n\} \backslash I$.
Note that under these assumptions, for all $x \in X$ we have $|g(x \cap I)| = |g(x) \cap I|$ and
$|g(x \cap J)| = |g(x) \cap J|$.  Moreover, we must have $g(I) = I$.

We will reach a contradiction with weak irreducibility along $X$  by showing that if $<i,j>\ \in A_x$ for some $x \in X$, then
we cannot have $i \in I$ and $j \in J$. Suppose that $<i,j>\ \in A_x$ for $x \in X$ with $i \in I$ and $j \in J$.
Then there exists $y$ such that either $x < x \cup \{i\} = y$, $j \in g(y)
\backslash g(x)$ or $y = x \backslash \{i\} < x$, $j \in g(x) \backslash g(y)$.  Wlog assume the former; the proof in the latter case is analogous.
Let $z = I \cup (x \cap J)$.  Then $g(z) \supseteq g(I) = I$, and $g(z) \supseteq g(y) \supseteq g(x) \cup \{j\}$.  It follows that
$|g(z)| \geq |I| + |g(x) \cap J| + 1$, since $j \notin g(x)$.  But this implies $|g(z)| > |z|$, contradicting the assumption of strong cooperativity.
\epr

For non-Boolean systems, the assumption of irreducibility in Theorem~\ref{itosc} or Corollary~\ref{shrtcyccor} cannot be replaced by the assumption of
strong semi-irreducibility.

\begin{Example}\label{almostex}
For every $n$ there exists a cooperative strongly semi-irreducible $4$-discrete system $(\Pi, g)$ of dimension $n$ that contains a periodic orbit
of length $d_{n,2}$.
\end{Example}

\bpr
Fix $n$, let $(\{0,1\}^n, f)$ be a cooperative Boolean system with a periodic orbit of length $d_{n,2}$, and let $\pi$ be a cyclic
permutation of $\{1, \ldots , n\}$.  Let $\Pi = \{0, 1, 2, 3\}^n$ and define a function $g: \Pi \rightarrow \Pi$ as follows.
Let $S = \{x \in \Pi: \ \min x = 0 \leq \max x < 3$, $M = \{x \in \Pi: \ 1 \leq \min x \leq \max x \leq 2\}$,
and $L = \{x \in \Pi: \ \max x = 3\}$.  For $x \in S$, let $g(x)_{\pi(i)} = 0$ whenever $x_i = 0$ and $g(x)_{\pi(i)} = 1$ whenever $x_i > 0$.
For $x \in M$, let $g(x)_j = 1 + f(x -1)$, and for $x \in L$ let  $g(x)_{\pi(i)} = 2$ whenever $0 < x_i < 3$,  $g(x)_{\pi(i)} = 0$
whenever $x_i = 0$, and $g(x)_{\pi(i)} = 3$ whenever
$x_i = 3$.

Note that the restriction  $g\rest M$ is isomorphic to $f$, hence the restriction of our system to $M$ is cooperative and has a periodic orbit
of length $d_{n,2}$.  It also follows immediately from the definitions that the restriction of our system to $S$ as well as its
restriction to $L$ are cooperative.  Moreover, consider $x \in S$, $y\in M$ and $z \in L$. Then $g(x) \leq g(y)$, and $x \leq z$ implies $g(x) \leq g(z)$.
Similarly, $y \leq z$ implies $g(y) \leq g(z)$. Since no
element of $S$ can sit above an element of $M$ or $L$, and no element of $M$ can sit above an element of $L$, strong cooperativity
of the whole system follows.

It remains to show that our system is strongly semi-irreducible.  It suffices to show that if $\pi(i) = j$, then $<i,j> \ \in A^*_x$
for all $x \in \Pi$.  Fix $x$ and $i, j$ with $j = \pi(i)$.  If $x_i = 3$, then $x \in L$ and $(g(x^{i-}))_{j} \leq 2 < 3 = (g(x))_j$.
If $x_i = 2$, then $x^{i+} \in L$ and $(g(x))_j \leq 2 < 3 = (g(x^{i+}))_j$.  Similarly, if
$x_i = 1$, then $x^{i-} \in S \cup L$ and $(g(x^{i-}))_j = 0 < 1 \leq (g(x))_j$.
Finally, if $x_i = 0$, then $x \in S \cup L$ and $(g(x))_j = 0 < 1 \leq (g(x^{i+}))_j$.
\epr

It turns out that for Boolean cooperative systems, strong irreducibility along an attractor $X$ all
by itself (without the assumption of strong cooperativity) puts tight bounds on the length of
this attractor.

\begin{Theorem}\label{sirshort}
Suppose $(\Pi, g)$ is an $n$-dimensional cooperative Boolean system, and let $X$ be an attractor such that $(\Pi, g)$ is strongly irreducible
along $X$.  Then $|X| \leq n$.
\end{Theorem}

\bpr
Let $(\Pi, g)$ and $X$ be as in the assumption, and let $A = \bigcap_{x \in X} A_x$. We will identify $x \in \Pi$ with subsets of $\{1, \ldots , n\}$. Note that if $<i, j> \ \in A$, then for all $x \in X$
we have $i \in x$ iff $j \in g(x)$.  By induction, if $j$ can be reached by a directed path in $G$ of length $r$, then
$i \in x$ iff $j \in g^r(x)$.  For each $1 \leq j \leq n$ define $\ell(j)$ as the length of the shortest directed
path in $G$ from $1$ to $j$.   Strong connectedness implies that $\ell(j)$ is well-defined and $1 \leq \ell(j) \leq n$ for all $j$.
Now let $x(0) \in X$ and note that $x(r\ell(1) + s)_1 = (g^{r\ell(1) + s}(x))_1 = (g^{s}(x))_1 = (x(s))_1$ for all nonnegative integers $r, s$.  More generally,
$x(r\ell(1) + s + \ell(j))_j = x(s)_1$ for all $r$, $s$~and $j$.  Thus $(x(r\ell(1)))_j = (x(-\ell(j)))_j = (x(0))_j$ for all
$j$, and hence $x(r\ell(1)) = x(0)$. The theorem follows.
\epr

Theorem~\ref{sirshort} does not generalize to arbitrary cooperative $p$-discrete systems.  In general, a $p$-discrete system may be
strongly irreducible even along an attractor of exponential length.

\begin{Example}\label{nopsirshortex}
For every $n$ there exists a cooperative $6$-discrete system $(\Pi, g)$ of dimension~$n$ that is strongly irreducible along a periodic orbit
of length $d_{n,2}$.
\end{Example}

\bpr
The construction is similar to Example~\ref{almostex}.
Fix $n$ and let $(\{0,1\}^n, f)$ be a cooperative Boolean system where $D = \{x \in \{0, 1\}^n: \ S(x) = \lfloor n/2\rfloor\}$ is a periodic orbit, and let $\pi$ be a cyclic
permutation of $\{1, \ldots , n\}$.  Let $\Pi = \{0, ... , 5\}^n$ and let $M = \{1 + 3x: \, x \in D\}$.
Define $h(1 + 3x) = 1 + 3f(x)$ for all $x \in D$.  For $y \in M$ and $i \in \{1, \ldots , n\}$, define
$h(y^{i+}) = h(y)^{\pi(i)+}$ and $h(y^{i-}) = h(y)^{\pi(i)-}$.  Let $A$ be the set of $y \in \Pi$ for which we have defined $h$.   Note that
for $x, y \in M$ and $i, j \in \{1, \ldots n\}$ we can have $x^{i-} \leq y^{j+}$ only if $x = y$.  Thus $h$ is cooperative on $A$.
By Lemma~\ref{smaleextension}, there exists a cooperative function
$g: \Pi \rightarrow \Pi$ so that $g \rest A = h$.  Clearly, $M$ is a periodic orbit of length $d_{n,2}$ of $g$, and it follows from
our definitions that $<i, \pi(i)> \ \in A_y$ for all $y \in M$.  Thus $(\Pi, g)$ is strongly irreducible along $M$.
\epr

The following example shows that assumption in Lemma~\ref{sirshort} cannot be weakened to irreducibility along the attractor.

\begin{Example}\label{irlong}
For every $n \geq 1$ there exists an $n$-dimensional cooperative Boolean system $(\Pi, g)$ and a periodic orbit $X$ of
length $d_{n,2}$ such that the system is irreducible along $X$.
\end{Example}

\bpr
The statement is trivially true for $n = 1$. So let  $n>1$, let $D$ be as in~(\ref{Ddef})
 and let $f: D \rightarrow D$ be such that the system $(D, f)$ is cyclic.
Embed $(D, f)$ into a cooperative Boolean system $(\Pi, g)$ as in Lemma~\ref{lemma trivial embedding}.  Consider $x \in D$.
Then $<i, j> \ \in A_x$ iff either $i \in x$ and $j \in f(x)$ or $i \notin x$ and $j \notin f(x)$. Let $i \in \{1, \ldots , n\}$.
We will show that every node in $\{1, \ldots , n\}$ can be reached from $i$ by a directed path in $G_x$ which implies
strong connectedness of $G_x$.
Assume $i \in x$; the proof in the case $i \notin x$ is dual. Then $<i, j> \ \in A_x$ for all $j \in f(x)$.
Now note that since $f(x) \neq  x$ and both $x, f(x) \in D$, there must  exist $j \in f(x) \backslash x$.
Then $<j, k > \ \in A_x$ for every $k \notin f(x)$, and it follows that every node in $\{1, \ldots , n\}$ can be reached from $i$ by a directed path
of length at most two.
\epr

It is not in general true for every periodic attractor of length $d_{n,2}$ in an $n$-dimensional cooperative Boolean system
that the system is irreducible along $X$.  We have the following example.

\begin{Example}\label{Germanex}
Fix $p >1$, let $n$ sufficiently large and let $D$ be as in~(\ref{Ddef}).
Then there exists a function $\gamma: D \rightarrow D$ such that whenever $(\Pi, g)$ is an $n$-dimensional cooperative
$p$-discrete system with $g \rest D = \gamma$, then $D$
is a periodic attractor of length $d_{n,p}$ such that $\bigcup_{x \in D} A^*_x$ contains all arcs $<i, j>$ with $i \neq j$.
In particular, the system $(\Pi, g)$ is weakly irreducible along $D$.  However, if $g$ is the Smale extension of $f$,
then $(\Pi,g)$ is not irreducible along $D$.
\end{Example}

\bpr
Let $<n>^2$ denote the set of all pairs $s=\ <i,j>$, for arbitrary $1\leq i<j\leq
n$. For any $s=\ <i,j>\ \in \ <n>^2$, let $\pi_s:D\to D$ be the result of swapping
the values of the components $i$ and $j$.  That is, for $d\in D$, let $(\pi_s\,
d)_i:=d_j$, $(\pi_s\, d)_j:=d_i$, and $(\pi_s\, d)_k:=d_k$ for all $k\not=i,j$.

Consider a function $a:\ <n>^2\to D$, and define $b:\ <n>^2\to D$ by $b(s):=\pi_s\,
a(s)$. Similarly, consider another function $a':\ <n>^2\to D$ and $b'(s):=\pi_s\,
a'(s)$. We prove the following result using the probabilistic method.

\begin{Lemma} \label{lemma injective}
For large enough $n$, the functions $a,a'$ above can be chosen in such a way
that $\mbox{Im } a \cup \mbox{Im } b \cup \mbox{Im } a' \cup \mbox{Im } b'$ has
exactly $4 \binom{n}{2}$ elements.
\end{Lemma}

\bpr
Let $T(s)$ denote the set $\{ a(s),b(s),a'(s),b'(s)\}$, for any $s\in<n>^2$.
The objective is to choose $a$ and $a'$ so that

\begin{equation}\label{counting 1}
\mbox{for every } s=\ <i,j>\ \in \ <n>^2:\ \  \norma{T(s)}=4,
\end{equation}
and so that
\begin{equation}\label{counting 2}
\mbox{for every $s,r\in \ <n>^2$, $s\not=r$: }T(s) \cap T(r) =\emptyset.
\end{equation}

For $s=\ <i,j>\ \in \ <n>^2$, define $D_s$ as the set of all $x\in D$ such that
$x_i>x_j$.  Let $q_n$ be the size of this set. Note that $q_n$ does not
depend of the particular values of $i<j$.

If $a(s),a'(s)\in D_s$, then an equivalent condition for~(\ref{counting 1}) is
that $a(s)\not=a'(s)$.  To see the sufficiency of this condition (the necessity
being obvious), note that if $a(s)\not=a'(s)$, then necessarily
$b(s)\not=b'(s)$ by construction, and that also $b(s),b'(s)\not\in D_s$, so
that $a(s)\not=b(s), b'(s)$; similarly $a'(s)\not=b(s), b'(s)$.

We will prove this lemma by making a simple use of the probabilistic method.
Suppose that for every fixed $s$, $a(s)$ is a discrete random variable that
takes any value in $D_s$ with equal probability.  Let $a'(s)$ be similarly
defined, and in such a way that all random values of both variables for any
$s,r\in \ <n>^2$ are independently distributed.  We show that the probability
that both (\ref{counting 1}) and (\ref{counting 2}) hold is greater than zero,
which implies the result.

For example, for any given $s$, we compute the probability that $a(s)=a'(s)$:
\[
P(a(s)=a'(s))=\sum_{x\in D_s} P(a(s)=x, a'(s)=x) =q_n
\frac{1}{q_n}\frac{1}{q_n}=\frac{1}{q_n}.\]
If $r\not=s$, we compute the probability that $a(s)=a(r)$:
\[
P(a(s)=a(r))=\sum_{x\in D_s\cap D_r} P(a(s)=x, a(r)=x) =\norma{D_s\cap D_r}
\frac{1}{q_n} \frac{1}{q_n} \leq \frac{1}{q_n}.
\]
This bound also holds for $a(s)=a'(r)$, and for any other combination of the
functions $a,a',b,b'$ (note that it is quite possible that, say, $a(s)=b(r)$
for $r\not=s$).  Thus we obtain for $s\not=r$ the bound:

\[
P(T(s) \cap T(r)\not=\emptyset) \leq 16 \frac{1}{q_n}.
\]

Now the event that either (\ref{counting 1}) or (\ref{counting 2}) fails can
be described as

\begin{equation} \label{P inequality}
\begin{array}{l}  \displaystyle
P(\mbox{either (\ref{counting 1}) or (\ref{counting 2}) fails}) \\
\leq   \displaystyle \sum_{s\in\ <n>^2} P(a(s)=a'(s))\  +  \sum_{r,s\in \ <n>^2, r\not=s} P(T(s)\cap T(r)\not=\emptyset) \\
\leq  \displaystyle \sum_{s\in \ <n>^2} 16 \frac{1}{q_n}\   +   \sum_{r,s\in \
<n>^2, r\not=s} 16 \frac{1}{q_n} = \binom{n}{2}(\binom{n}{2} + 1) 16 \frac{1}{q_n}
\end{array}
\end{equation}

We give a lower bound for $q_n$ as follows.  Consider the set $D'$ of vectors
$x\in \{0,\ldots p-1\}^{n-2}$ such that $S(x)=\flo{(n-2)(p-1)/2}$, which has
cardinality $d_{n-2,p}$.  A straightforward injection of $D'$ into
$D_{(n-1,n)}$ is given by the function $[x_1,\ldots, x_{n-2}]\to [x_1,\ldots,
x_{n-2},p-1,0]$ (note $\flo{(n-2)(p-1)/2}+(p-1)=\flo{\frac{n(p-1)}{2}}$).

This
implies that
\[
q_n\geq d_{n-2,p} \geq \frac{p^{n-3}}{n},
\]
by Lemma~\ref{lemma dnp lower bound}.  Since $ \binom{n}{2}(\binom{n}{2} + 1)
\frac{16}{q_n}\leq \binom{n}{2}(\binom{n}{2} + 1) \frac{16 n}{p^{n-3}}<1$  for large enough
$n$, the statement follows from (\ref{P inequality}).
\epr

Now consider functions $a,a':\ <n>^2\to D$ as in the statement of Lemma~\ref{lemma
injective}, and their associated functions $b,b'$ such that the conclusion of
this Lemma is satisfied.  In particular, for every $s= \ <i,j>$,
$a(s)_i\not=a(s)_j$ must hold (else $a(s)=b(s)$).  By swapping the values of
$a(s)$ and $b(s)$ if necessary, we can assume without loss of generality that
$a(s)_i>a(s)_j$, $b(s)_i< b(s)_j$ for every $s$.  Repeat the same procedure
with $a',b'$.

Define $c,c':\ <n>^2\to \Pi$ by $c(s)=min(a(s),b(s))$, $c'(s)=min(a'(s),b'(s))$.
In particular $c(s)<a(s),\ c(s)<b(s)$ for every $s$, and similarly with $c'$.

Define the values of $\gamma$ in such a way that $\gamma(a(s))=b'(s)$,
$\gamma(b(s))=a'(s)$ for every~$s$.  This is possible by Lemma~\ref{lemma
injective}, since the states $a(s)$, $b(s)$ are all different from each other.
Since we are still free to choose the values $\gamma(a'(s))$, $\gamma(b'(s))$
for every $s$ (again by Lemma~\ref{lemma injective}), we can define the values
of $\gamma$ in such a way that $\gamma$ generates a single orbit of period
$\norma{D}$.

Let $g:\Pi\to \Pi$ be a cooperative extension of $\gamma$.  For any fixed $s= \ <i,j>$, since $c(s)\leq b(s)$, it follows by
cooperativity $g(c(s))\leq g(b(s))$.  But then

\[
g(c(s))_j\leq g(b(s))_j=a'(s)_j < b'(s)_j=g(a(s))_j.
\]
Since $c(s)_k=a(s)_k$ for all $k\not=i$, and $c(s)_i<a(s)_i$, the equation
above shows that $<i,j> \ \in A^*_{a(s)}$.  Similarly, using $c(s)\leq a(s)$,
\[
g(c(s))_i\leq g(a(s))_i=b'(s)_i < a'(s)_i=g(b(s))_i,
\]
which in the same way implies $<j,i> \ \in A^*_{a(s)}$.  This shows that $\bigcup_{x \in D} A^*_x$ contains all arcs $<i, j>$ with $i \neq j$;
in particular, the system $(\Pi, g)$ is weakly irreducible along $D$.

Finally, assume that $g$ is the Smale extension of $\gamma$ as defined in Section~\ref{smaleapprox}.
Note that for every $s$, the only elements $d \in D$ with $c(s) \leq d$ are $a(s)$ and $b(s)$.  Thus
$g(c(s)) = \inf \{a'(s), b'(s)\}$, and our choice of these values implies that
$(g(c(<i,j>)))_k = (g(a(<i,j>)))_k$ for all $k \neq i, j$.  Similarly, if $e(s) = \max \{a(s), b(s)\}$, then
$(g(c(<i,j>)))_k = (g(a(<i,j>)))_k$ for all $k \neq i, j$. Thus $A_{a(<i,j>)}$ contains no arcs from the set
$\{i, j\}$ into its complement, and thus $G_{a(<i,j>)}$ cannot be strongly connected when $n > 2$.
It follows that $(\Pi, g)$ cannot be irreducible along $D$.
\epr

\textbf{Acknowledgements:}  The authors would like to thank Eduardo Sontag for
many helpful comments and suggestions.  We also want to thank Neil Falkner for
providing a reference to the proof of Proposition~\ref{proposition local
limit}, as well as Akos Seress, David Terman, and Tom Wolf for helping us to locate
the reference to the maximal order of permutation groups. 

%\bibliography{DissertationRefs}
%\bibliographystyle{plain}

\end{document}